\input lanlmac
\input epsf.tex
\input mssymb.tex
\font\sc=cmcsc10
\input threearches.defs
\overfullrule=0pt

\newcount\figno
\figno=1
\def\fig#1#2#3{%
\xdef#1{\the\figno}%
\writedef{#1\leftbracket \the\figno}%
\nobreak%
\par\begingroup\parindent=0pt\leftskip=1cm\rightskip=1cm\parindent=0pt%
\baselineskip=11pt%
\midinsert%
\centerline{#3}%
\vskip 12pt%
{\bf Fig.\ \the\figno:} #2\par%
\endinsert\endgroup\par%
\goodbreak%
\global\advance\figno by1%
}
\long\def\rem#1{{\sl [#1]}}
\def\e#1{{\rm e}^{#1}}
\def\pre#1{{\tt
#1}}
\def\d{{\rm d}}

\def\pf{{\rm Pf}}

%
%
\newcount\propno
\propno=1
\def\prop#1#2\par{\xdef#1{\the\propno}%
\medbreak\noindent{\sc Proposition \the\propno.\enspace}{\sl #2}\medbreak%
\global\advance\propno by1}
\def\thmn#1#2\par{
\medbreak\noindent{\sc Theorem #1.\enspace}{\sl #2\par}\medbreak%
}
\newcount\thmno
\thmno=1
\def\thm#1#2\par{\xdef#1{\the\thmno}%
\thmn{\the\thmno}{#2}%
\global\advance\thmno by1}
\newcount\conjno
\conjno=1
\def\conj#1#2\par{\xdef#1{\the\conjno}%
\medbreak\noindent{\sc Conjecture \the\conjno.\enspace}{\sl #2}\medbreak%
\global\advance\conjno by1}
\newcount\lemno
\lemno=1
\def\lemma#1#2\par{\xdef#1{\the\lemno}%
\medbreak\noindent{\sc Lemma \the\lemno.\enspace}{\sl #2}\medbreak%
\global\advance\lemno by1}
\def\corol#1\par{%
\medbreak\noindent{\sc Corollary.\enspace}{\sl #1}\medbreak}
\def\example#1\par{%
\medbreak\noindent{\sc Example:\enspace}#1\medbreak}
\def\qed{\nobreak\hfill\vbox{\hrule height.4pt%
\hbox{\vrule width.4pt height3pt \kern3pt\vrule width.4pt}\hrule height.4pt}\medskip\goodbreak}
\lref\RS{A.V. Razumov and Yu.G. Stroganov, 
{\sl Combinatorial nature
of ground state vector of $O(1)$ loop model},
{\it Theor. Math. Phys.} 
{\bf 138} (2004) 333--337; {\it Teor. Mat. Fiz.} 138 (2004) 395--400, \pre{math.CO/0104216}.}
\lref\BdGN{M.T. Batchelor, J. de Gier and B. Nienhuis,
{\sl The quantum symmetric XXZ chain at $\Delta=-1/2$, alternating sign matrices and 
plane partitions},
{\it J. Phys.} A34 (2001) L265--L270,
\pre{cond-mat/0101385}.}
\lref\PRdG{P. A. Pearce, V. Rittenberg and J. de Gier, 
{\sl Critical Q=1 Potts Model and Temperley--Lieb Stochastic Processes},
\pre{cond-mat/0108051}.}
\lref\RSb{A.V. Razumov and Yu.G. Stroganov, 
{\sl $O(1)$ loop model with different boundary conditions and symmetry classes of alternating-sign matrices},
{\it Theor. Math. Phys.} 
{\bf 142} (2005) 237--243; {\it Teor. Mat. Fiz.} 142 (2005) 284--292,
\pre{cond-mat/0108103}.}
\lref\PRdGN{P. A. Pearce, V. Rittenberg, J. de Gier and B. Nienhuis,
{\sl Temperley--Lieb Stochastic Processes}, 
{\it J. Phys. A} {\bf 35 } (2002) L661-L668, \pre{math-ph/0209017}.}
\lref\dG{J.~de~Gier, {\sl Loops, matchings and alternating-sign matrices},
\pre{math.CO/0211285}.}
\lref\MNosc{S. Mitra and B. Nienhuis, {\sl 
Osculating random walks on cylinders}, in
{\it Discrete random walks}, 
DRW'03, C. Banderier and
C. Krattenthaler edrs, Discrete Mathematics and Computer Science
Proceedings AC (2003) 259-264, \pre{math-ph/0312036}.} 
\lref\MNdGB{S. Mitra, B. Nienhuis, J. de Gier and M.T. Batchelor,
{\sl Exact expressions for correlations in the ground state 
of the dense $O(1)$ loop model}, 
JSTAT (2004) P09010,
\pre{cond-math/0401245}.}
\lref\DF{P.~Di Francesco, {\sl 
 A refined Razumov--Stroganov conjecture} I: 
     J. Stat. Mech. P08009 (2004), \pre{cond-mat/0407477}; II: 
     J. Stat. Mech. P11004 (2004), \pre{cond-mat/0409576}.}
\lref\LGV{B. Lindstr\"om, {\it On the vector representations of
induced matroids}, Bull. London Math. Soc. {\bf 5} (1973)
85--90\semi
I. M. Gessel and X. Viennot, {\it Binomial determinants, paths and
hook formulae}, Adv. Math. { \bf 58} (1985) 300--321. }
\lref\DFZJ{P.~Di Francesco and P.~Zinn-Justin, {\sl Around the Razumov--Stroganov conjecture:
proof of a multi-parameter sum rule}, {\it E. J. Combi.} 12 (1) (2005), R6,
\pre{math-ph/0410061}.}
\lref\Pas{V.~Pasquier, {\sl Quantum incompressibility and Razumov Stroganov type conjectures},
\pre{cond-mat/0506075}.}
\lref\DFZJc{P.~Di Francesco and P.~Zinn-Justin, 
{\sl Quantum Knizhnik--Zamolodchikov equation, generalized Razumov--Stroganov sum rules 
and extended Joseph polynomials}, 
{\it J. Phys. A} 38 (2005) L815--L822, \pre{math-ph/0508059}.}
\lref\DFZJZ{P.~Di~Francesco, P.~Zinn-Justin and J.-B.~Zuber,
{\sl A Bijection between classes of Fully Packed Loops and Plane Partitions},
{\it E. J. Combi.} 11(1) (2004), R64,
\pre{math.CO/0311220}.}
\lref\DFZ{P.~Di~Francesco and J.-B.~Zuber, {\sl On FPL 
configurations with four sets of nested arches}, 
JSTAT (2004) P06005, \pre{cond-mat/0403268}.}
\lref\DFZJZb{P.~Di~Francesco, P.~Zinn-Justin and J.-B.~Zuber,
{\sl Determinant Formulae for some Tiling Problems and Application to Fully Packed Loops}, 
{\it Annales de l'Institut Fourier} 55 (6) (2005), 2025--2050, \pre{math-ph/0410002}.} 
\lref\KOR{V. Korepin, {\sl Calculation of norms of Bethe wave functions},
{\it Comm. Math. Phys.} {\bf 86} (1982) 391-418.}
\lref\Kratt{F. Caselli and C.~Krattenthaler, {\sl Proof of two conjectures of
Zuber on fully packed loop configurations}, {J. Combin. Theory
Ser.} {\bf A 108} (2004), 123--146, \pre{math.CO/0312217}.  }
\Title{}
{
\vbox{\centerline{Proof of Razumov--Stroganov conjecture}
\medskip
\centerline{for some infinite families of link patterns}}
}
\bigskip\bigskip
\centerline{P. Zinn-Justin \footnote{${}^\star$}
{Laboratoire de Physique Th\'eorique et Mod\`eles Statistiques, UMR 8626 du CNRS,
Universit\'e Paris-Sud, B\^atiment 100,  F-91405 Orsay Cedex, France}}
\vskip0.5cm
\noindent
We prove the Razumov--Stroganov conjecture relating ground state of the $O(1)$ loop model and counting of
Fully Packed Loops in the case of certain types of link patterns. The main focus is on link
patterns with three series of nested arches, for which we use as key ingredient of the proof a
generalization of the Mac Mahon formula for the number of plane partitions which includes three
series of parameters.

\bigskip

\long\def\rem#1{}

\Date{04/2006}
%
%
%
\newsec{Introduction}
The Razumov--Stroganov (RS) conjecture \RS\ relates the components of the ground state of the
$O(1)$ loop model, which are indexed by link patterns (pairing of points
on a circle) to the numbers of Fully Packed Loop configurations (FPL) on a square grid
with a connectivity of external vertices given by the link pattern. Despite considerable activity around
this conjecture \refs{\PRdG,\RSb,\PRdGN,\dG,\MNosc,\MNdGB,\DF},
it has not been proved yet. It is the author's belief, however, that the work \DFZJ\ 
was a significant step in this direction: in it, an inhomogeneous loop model was introduced in order
to make the ground state a polynomial of the inhomogeneities (spectral parameters). This
way, a corollary of the RS conjecture (which was already formulated in \BdGN), namely that the properly
normalized sum of all
components of the loop model ground state equals the total number of FPL, also known as the number of Alternating
Sign Matrices, was proved.

The present work tries to demonstrate in a very simple setting how the methods of \DFZJ\ could help to prove
the RS conjecture by considering a special subset of possible link patterns, namely
those with few ``little arches'' (arches connecting neighbors). Consider the link patterns of Fig.~\lpthree.
They are made of three sets of $a$, $b$, $c$ nested arches.
Here $a$, $b$, $c$ are arbitrary integers such that $a+b+c=n$ where the size of the system is $2n$.
The model depends on $2n$ complex numbers $\alpha_1,\ldots,\alpha_{b+c},\beta_1,\ldots,\beta_{a+c},\gamma_1,\ldots,
\gamma_{a+b}$ which are, up to multiplication by a power of $q$ (as will be explained below), the
spectral parameters of the model.
Note there are several good reasons to restrict oneself to such link patterns, a particularly obvious
one being that we know the corresponding number of FPL: it was computed in \DFZJZ\ -- and happens
to be equal to the number of Plane Partitions in a hexagon of shape $a\times b\times c$!
\fig\lpthree{Link pattern with 3 sets of nested arches.}{\epsfbox{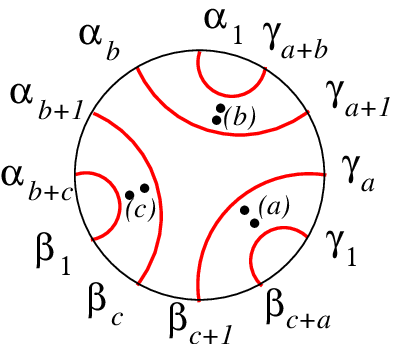}}

The reader is referred to \DFZJ\ for details concerning the $O(1)$ loop model. We briefly recall its
definition here, if only to fix notations. Let $n$ be an integer, and
consider a complex vector space equipped with a basis indexed
by {\it non-crossing link patterns}: the latter are by definition pairings of $2n$ points on a circle, 
in such a way that the pairings can be represented by non-crossing edges inside the circle. On this
space we define the action of a linear operator, the so-called transfer matrix $T_n(t\vert z_1,\ldots,z_{2n})$, 
which depends on complex parameters $t,z_1,\ldots,z_{2n}$, by the following graphical description:
\eqn\tm{
T_n(t\vert z_1,\ldots,z_{2n})=\prod_{i=1}^{2n} \left(
{q\,z_i-q^{-1}t\over q\,t-q^{-1}z_i}
\vcenter{\hbox{\epsfbox{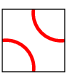}}}
+
{z_i-t\over q\,t-q^{-1}z_i}
\vcenter{\hbox{\epsfbox{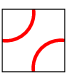}}}
\right)
}
where $q=\e{2i\pi/3}$, and the symbolic product over $i$ means that the plaquette of index $i$ 
should be inserted at vertex $i$. 
The result is that, starting from a given link pattern, the action of $T_n$ produces
a new link pattern by adding a circular strip of plaquettes and removing any closed loops thus
created.
The coefficients in Eq.~\tm, in the range of parameters where they are real, are simply
probabilities of inserting the corresponding plaquettes, parametrized in a convenient way in terms of the
$z_i/t$. As a consequence of the Yang--Baxter equation, $[T_n(t),T_n(t')]=0$ with all other parameters
$z_i$ fixed.

Since $T_n$ is a stochastic matrix, it has the obvious left eigenvector $(1,\ldots,1)$ with eigenvalue $1$.
Therefore, it also has a corresponding right eigenvector, which is unique for generic values of the $z_i$:
\eqn\pf{
T_n(t\vert z_1,\ldots,z_{2n}) \Psi_n(z_1,\ldots,z_{2n})=\Psi_n(z_1,\ldots,z_{2n})
}
One can normalize $\Psi_n$ in such a way that its components $\Psi_{n,\pi}$ in the basis of link patterns $\pi$
are {\it coprime polynomials} of the $z_i$. One still has an arbitrary numerical constant in the normalization
of $\Psi_n$. Consider now the homogeneous limit when all $z_i$ equal $1$. If one chooses this constant
so that the smallest entry is $1$, then the $\Psi_{n,\pi}(1,\ldots,1)$ are the subject of various conjectures,
including the remarkable Razumov--Stroganov conjecture already mentioned above, that identifies them
with a certain FPL enumeration problem.
In what follows we shall choose another numerical normalization which is more convenient for intermediate
calculations.

In Sect.~2, we derive the main formula for the entries of the ground state of the $O(1)$ loop model corresponding
to link patterns with three sets of nested arches (Fig.~\lpthree). 
In Sect.~3 we establish the connection with plane partitions 
(or dimers). Sect.~4 discusses the (partial or total)
homogeneous limit. Sect.~5 briefly describes the extension to 
more general link patterns
for which the corresponding enumeration of FPL is known. Sect.~6 concludes.

\newsec{Recurrence relations and their solution}
In \DFZJ, a certain number of relations were shown to be satisfied
by $\Psi_n$, the ground state eigenvector of the $O(1)$ loop model
of size $2n$. We shall need the following three facts:

\thmn{4 of \DFZJ}The components of $\Psi_n$ are homogeneous
polynomials of total degree $n(n-1)$, and of partial degree at most $n-1$ in each variable $z_i$.

\thmn{1 of \DFZJ}The entries $\Psi_{n,\pi}$ of the groundstate
eigenvector satisfy:
\eqn\recone{\Psi_{n,\pi}(z_1,\ldots,z_{2n})=
\Bigg(\prod_{s\in E_\pi} \prod_{\scriptstyle i,j\in s\atop
\scriptstyle i<j} (q z_i-q^{-1}z_j)\Bigg)
\Phi_{n,\pi}(z_1,\ldots,z_{2n})}
where $\Phi_{n,\pi}$ is a polynomial which is symmetric in the set of variables
$\{ z_i, i\in s\}$ for each $s\in E_\pi$, and $E_\pi$ is the partition of
$\{1,\ldots,2n\}$ into maximal sequences of consecutive points not connected to each other by arches of $\pi$.

\thmn{3 of \DFZJ}If two neighboring parameters $z_i$ and $z_{i+1}$ are
such that $z_{i+1}=q^2 z_i$, then either of the two following situations occur for the components
$\Psi_{n,\pi}$:\hfill\break
(i) the pattern $\pi$ has no arch joining $i$ to $i+1$,
in which case according to Theorem 1,
\eqn\vanish{\Psi_{n,\pi}(z_1,\ldots,z_i,z_{i+1}=q^2 z_i,\ldots,z_{2n})=0\ ;}
(ii) the pattern $\pi$ has a little arch joining $i$ to $i+1$, 
in which case
\eqnn\recure
$$\eqalignno{\Psi_{n,\pi}(z_1,&\ldots,z_i,z_{i+1}=q^2 z_i,\ldots,z_{2n})=\cr
&
\left(\prod_{\scriptstyle k=1\atop\scriptstyle k\ne i,i+1}^{2n} (q\,z_i-z_k)\right)
\ \Psi_{n-1,\pi'}(z_1,\ldots,z_{i-1},z_{i+2},\ldots,z_{2n})&\recure}$$
where $\pi'$ is the link pattern $\pi$ with the little arch $i$, $i+1$
removed.

In what follows we shall concentrate on components corresponding to link patterns 
with three sets of nested arches
of size $a,b,c$, which we shall denote by $\Psi_{a,b,c}$. The spectral parameters are relabelled
as $z_i=\alpha_i,q\,\beta_i,q^2\gamma_i$ according to the pattern of Fig.~\lpthree.
Thanks to Theorem 1, we can write
\eqn\facto{
\Psi_{a,b,c}=
\prod_{1\le i<j\le b+c}\!\!(q\,\alpha_i-q^{-1}\alpha_j)
\prod_{1\le i<j\le a+c}\!\!(q\,\beta_i-q^{-1}\beta_j)
\prod_{1\le i<j\le a+b}\!\!(q\,\gamma_i-q^{-1}\gamma_j)
\ \Phi_{a,b,c}
}
where the arguments $\{z_i\}=\{\alpha_i, q\beta_i, q^2\gamma_i\}$ have been suppressed for brevity.
According to Theorems 1 and 4,
$\Phi_{a,b,c}$ is a polynomial of total degree $ab+bc+ca$, 
and a symmetric polynomial of the $\{\alpha_i\}$ of degree at most $a$ in each, of
the $\{\beta_i\}$ of degree at most $b$ in each,
and of the $\{\gamma_i\}$ of degree at most $c$ in each.

We now rewrite Eq.~\recure\ in the case when $z_i$ is the last parameter $\alpha$ and
$z_{i+1}$ is the first parameter $\beta$, in terms of $\Phi_{a,b,c}$. Since the latter is a symmetric
function it is actually irrelevant which $\alpha$ and which $\beta$ are singled out, and the result is:
\eqn\newrecur{
\Phi_{a,b,c}|_{\beta_j=\alpha_i}=\prod_{k=1}^{a+b} (\alpha_i-\gamma_k)\ \Phi_{a,b,c-1}
}
where the parameters of $\Phi_{a,b,c-1}$ are the same as those of $\Phi_{a,b,c}$, except
$\alpha_j$ and $\beta_i$ are removed. Since $\Phi_{a,b,c}$ is of degree $a$ in $\alpha_i$,
the equations \newrecur\ with $j=1,\ldots,a+c$ and fixed $i$ determine entirely $\Phi_{a,b,c}$ as long
as $c\ge 1$. They form a very simple recurrence relation which is supplemented
by the initial
condition $\Phi_{a,b,0}$: this corresponds to the so-called ``base link pattern'' (Fig.~\base), which is entirely
factorized by Theorem 1:
\eqn\startrecur{
\Phi_{a,b,0}=\prod_{i=1}^b\prod_{j=1}^a (\alpha_i-\beta_j)
}
\fig\base{Base link pattern.}{\epsfbox{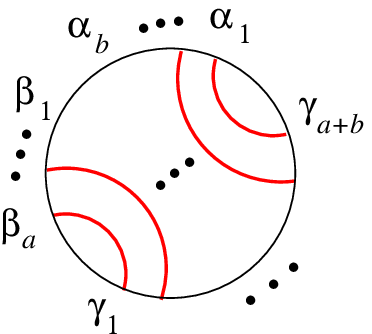}}

We now claim that the following Ansatz solves the recurrence relation:
\eqn\theansatz{
\Phi_{a,b,c}=\sum_{\scriptstyle I\subset\{1,\ldots,b+c\}\atop \# I=c}
{\prod_{i\not\in I}\prod_{j=1}^{a+c} (\alpha_i-\beta_j)
\prod_{i\in I}\prod_{k=1}^{a+b} (\alpha_i-\gamma_k)
\over
\prod_{i\in I}\prod_{j\not\in I} (\alpha_i-\alpha_j)
}
}

The case $c=0$ forces $I=\emptyset$ and we recover immediately Eq.~\startrecur. Next,
notice that the expression of Eq.~\theansatz\ is symmetric in all three sets of variables:
it is obvious for the $\{\beta_i\}$ and the $\{\gamma_i\}$; it is also clear for the $\{\alpha_i\}$ since the
summation over all possible subsets of cardinality $c$ is invariant by permutation of $\{1,\ldots,b+c\}$.
We can therefore choose one $\alpha_i$ and one $\beta_j$ and set them equal, say $\beta_{a+c}=\alpha_{b+c}$. 
This forces $b+c\in I$ in Eq.~\theansatz. Define $I'=I-\{b+c\}$. The summation
over $I$ can then be replaced with the summation over $I'\subset\{1,\ldots,b+c-1\}$, and it is easy
to check that cancellations in numerator and denominator reproduce Eq.~\theansatz\ with $c\to c-1$
and the parameters $\alpha_{b+c}$, $\beta_{a+c}$ removed. This concludes the recurrence.

Note immediately the symmetry between the sets of variables $\{\beta_i\}$ and $\{\gamma_i\}$ in 
Eq.~\theansatz: indeed, replacing $I$ with its complement $\bar I$ exchanges their roles (as well
as $b$ and $c$). However, the $\{\alpha_i\}$ seem to play a different role. We shall now produce
an equivalent expression which restores the symmetry $\alpha\leftrightarrow\beta$, at the expense
of breaking the symmetry $\beta\leftrightarrow\gamma$:
\eqnn\integrep
$$\eqalignno{
\Phi_{a,b,c}={1\over c!}\prod_{i=1}^{b+c}\prod_{j=1}^{a+c} (\alpha_i-\beta_j)
\oint&\cdots\oint {\d z_1\over 2\pi i}\cdots{\d z_c\over 2\pi i}
\prod_{1\le i< j\le c}(z_i-z_j)^2 \cr
&{\prod_{\ell=1}^c \prod_{k=1}^{a+b}(z_\ell-\gamma_k)
\over
\prod_{\ell=1}^c \prod_{j=1}^{a+c}(z_\ell-\beta_j)
\prod_{\ell=1}^c \prod_{i=1}^{b+c}(z_\ell-\alpha_i)
}&\integrep}$$
The $c$ contour integrals should be defined in such a way as to encircle (counterclockwise) all the poles
$\alpha_i$ (but none of the $\beta_i$). One goes back to Eq.~\theansatz\ by
applying the Cauchy formula. Each $z_i$ must be evaluated at a certain $\alpha_{I_i}$ with $1\le I_i\le b+c$;
furthermore, the factors $\prod (z_i-z_j)$ force the $I_i$ to be distinct, and we reproduce
after various cancellations the summation over $I=\{I_1,\ldots,I_{b+c}\}$ of Eq.~\theansatz.

The formula \integrep\ is of the form of a matrix integral: the contour integral makes it essentially
similar to the unitary matrix integral. This analogy will be pursued below. For now, we use a standard trick
in random matrix theory, which is to introduce the Vandermonde determinant
$\Delta(z_i)=\prod_{i<j}(z_i-z_j)=\det(z_i^{j-1})_{1\le i\le c}$, and then to note that
$\det(P_i(z_j))=\Delta(z_i)\det P$ where the $P_i$, $1\le i\le c$, are arbitrary polynomials
of degree less than $c$ and $P$ is the $c\times c$ matrix of coefficients of the $P_i$. 
In Eq.~\integrep\ we have a squared Vandermonde determinant, so
we can introduce another similar set of polynomials $Q_i$.
Moving the determinants out of the integrals, we find:
\eqn\detrep{
\Phi_{a,b,c}={\prod_{i=1}^{b+c}\prod_{j=1}^{a+c} (\alpha_i-\beta_j)\over \det P\det Q}
\det\left[
\oint {\d z\over 2\pi i}\, P_\ell(z)Q_m(z) 
{\prod_{k=1}^{a+b}(z-\gamma_k)
\over
\prod_{j=1}^{a+c}(z-\beta_j)
\prod_{i=1}^{b+c}(z-\alpha_i)
}\right]_{1\le \ell,m\le c}
}
In what follows, we shall be naturally led to a choice of polynomials $P$ and $Q$.

\newsec{Connection with plane partitions}
We now introduce a model of weighted Plane Partitions -- in more physical terms, it is a model
of dimers on the hexagonal lattice, but we shall prefer the language of Plane Partitions in what follows. 
Configurations are defined as tilings with lozenges of a hexagon of size $a\times b\times c$.
Lozenges come in three orientations since they are made of two adjacent equilateral triangles of
a regular triangular lattice.
The model comes with three series of parameters $\alpha_i$, $\beta_j$, $\gamma_k$
living on the lines of the underlying medial Kagome lattice, see Fig.~\ppfig~(i).
To each lozenge of the
plane partition (or equivalently to each dimer)
is associated a local Boltzmann weight 
$\alpha_i-\beta_j$, $\gamma_k-\beta_j$, $\alpha_i-\gamma_k$ given by the difference of the
parameters of the lines crossing at its center.
Note that there are exactly $ab$, $bc$, $ca$ lozenges of each orientation.
\fig\ppfig{(i) Plane partition and its parameters, $a=2$, $b=4$, $c=3$. (ii) Associated non-intersecting 
paths.}{(i) \epsfbox{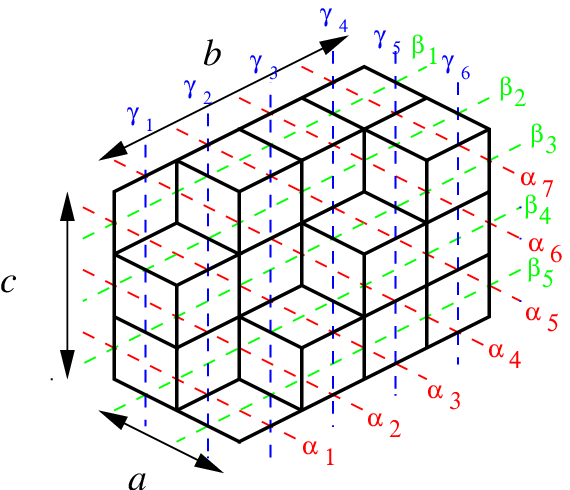}\qquad(ii) \epsfbox{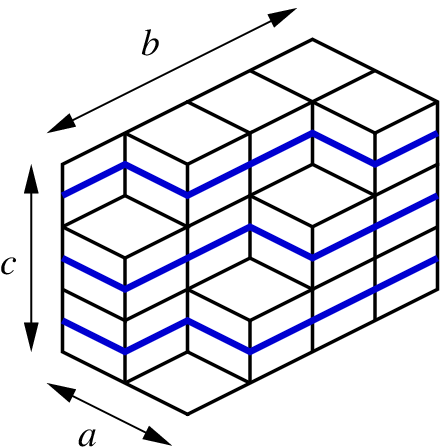}}

We wish to compute the partition function $Z_{a,b,c}$, 
i.e.\ the sum over all configurations of the product of local Boltzmann weights.
In order to do so, it is convenient to use yet another representation in terms of non-intersecting paths (NIPs).
To each plane partition one can associate $c$ paths originating from one of the sides of length $c$ and ending at the other, see
Fig.~\ppfig~(ii), which simply follow two types of tiles out of the three. 

One can replace
the local Boltzmann weights of the tiles with a local probability for the path to go left or right,
by factoring out all the possible weights of the third type of tile:
\eqn\partfun{
Z_{a,b,c}=
\prod_{
\scriptstyle 1\le i\le b+c, 1\le j\le a+c
\atop
\scriptstyle i+j>c,\ i+j\le a+b+c}
\!\!\!\!\!\!
(\alpha_i-\beta_j)\ \ F_{a,b;c}
}
and then introducing an inverse weight when the tiles of the third type are absent (i.e.\ where the paths
are):
\eqn\partfunb{
F_{a,b;c}=
\sum_{\rm NIPs}\,\prod_{{\rm edge}\in{\rm path}} \cases{
\displaystyle{\alpha_i-\gamma_k\over\alpha_i-\beta_j}&edge at the crossing of $(\alpha_i,\gamma_k)$\cr
\displaystyle{\gamma_k-\beta_j\over\alpha_i-\beta_j}&edge at the crossing of $(\beta_j,\gamma_k)$\cr
}
}
where the two orientations of the edges (or of the underlying dimers, or lozenges) determine which weight to use, and
the third coordinate is given by $i+j=k+c$. 
These NIPs move exactly $a$ steps in one direction and $b$ steps in the other.

NIPs are free fermions, and therefore their propagator $F_{a,b;c}$ is a determinant of one-particle propagators:
\eqn\lgv{
F_{a,b;c}=\det\left[\ell\to m\ {\rm on}\ (a,b,c)\right]_{1\le \ell,m\le c}
}
where $\ell\to m$ means the probability (for a single path) to go from position $\ell$ on one side of length $c$ to position $m$ on the other.
This is also known as the Lindstr\"om--Gessel--Viennot formula \LGV. 
The one-particle propagator is nothing but the case $c=1$, with appropriately shifted $a$, $b$, $\alpha_i$,
$\beta_j$: 
\eqn\oneparta{
[\ell\to m\ {\rm on}\ (a,b,c)]=F_{a+\ell-m,b+m-\ell,1}
(\alpha_\ell,\ldots,\alpha_{b+m};
\beta_{c+1-\ell},\ldots,\beta_{a+c+1-m}; \gamma_1,\ldots,\gamma_{a+b})
)
}

We are finally led to a simple problem of computing the weighted enumeration of a single path. 
The following formula holds:
\eqn\onepartb{
F_{a,b;1}(\alpha_1,\ldots,\alpha_{b+1};\beta_1,\ldots,\beta_{a+1};\gamma_1,\ldots,\gamma_{a+b})
=(\alpha_{b+1}-\beta_{a+1})\oint {\d z\over 2\pi i} {\prod_{k=1}^{a+b}(z-\gamma_k)\over\prod_{i=1}^{b+1}(z-\alpha_i)
\prod_{j=1}^{a+1}(z-\beta_j)}
}
where once again the contour integral encircles clockwise the $\alpha_i$ but not the $\beta_i$.
This can be proved by noting that $F_{a,b;1}$ satisfies the following simple recurrence formula:
\eqn\reconepart{
F_{a,b;1}={\alpha_{b+1}-\gamma_{a+b}\over \alpha_{b+1}-\beta_{a}}F_{a-1,b;1}
+{\gamma_{a+b}-\beta_{a+1}\over \alpha_{b}-\beta_{a+1}} F_{a,b-1;1}
}
and by using 
$z-\gamma_{a+b}=
{\alpha_{b+1}-\gamma_{a+b}\over \alpha_{b+1}-\beta_{a+1}}(z-\beta_{a+1})
+{\gamma_{a+b}-\beta_{a+1}\over \alpha_{b+1}-\beta_{a+1}}(z-\alpha_{b+1})$.

Putting together Eqs.~\partfun--\onepartb, one obtains
\eqn\detrepb{
Z_{a,b,c}=
\prod_{
\scriptstyle 1\le i\le b+c, 1\le j\le a+c
\atop
\scriptstyle i+j>c,\ i+j\le a+b+c+1}
\!\!\!\!\!\!
(\alpha_i-\beta_j) 
\ \ \det\left[
\oint {\d z\over 2\pi i}
{\prod_{k=1}^{a+b}(z-\gamma_k)
\over
\prod_{i=\ell}^{b+m}(z-\alpha_i)
\prod_{j=c+1-\ell}^{a+c+1-m}(z-\beta_j)
}\right]_{1\le \ell,m\le c}
}
To connect with Eq.~\detrep, set
\eqnn\polys
$$\eqalignno{
P_\ell(z)&=(z-\alpha_1)\cdots(z-\alpha_{\ell-1})\,(z-\beta_1)\cdots(z-\beta_{c-\ell})\cr
Q_m(z)&=(z-\alpha_{b+m+1})\cdots(z-\alpha_{b+c})\,(z-\beta_{a+c+2-m})\cdots(z-\beta_{b+c})&\polys\cr
}$$
By factor exhaustion one finds immediately that $\det P=\prod_{i+j\le c}(\alpha_i-\beta_j)$
and $\det Q=\prod_{i+j\ge a+b+c+2}(\alpha_i-\beta_j)$. Plugging this into
Eq.~\detrep\ reproduces exactly Eq.~\detrepb.

We conclude that $Z_{a,b,c}=\Phi_{a,b,c}$. We have thus obtained a direct, 
exact relation between the components of the inhomogeneous $O(1)$ loop model corresponding
to three sets of nested arches $(a,b,c)$ and the partition function of weighted plane partitions on a hexagon $a\times b\times c$.

\newsec{Relation to unitary matrix integrals and homogeneous limit}
In order to prepare the ground for the homogeneous limit, 
one can now choose some of the variables to be equal: $\alpha_i=\alpha$, $\beta_j=\beta$.
An important property of the matrix integral-like expression \integrep\ is that it is preserved by homographic
transformations. We define $w={z-\alpha\over z-\beta}$ to send
$(\alpha,\beta)$ to $(0,\infty)$ and obtain
\eqnn\integrepc
$$\eqalignno{
\Phi_{a,b,c}&={C\over c!}
\oint\cdots\oint {\d w_1\over 2\pi i}\ldots{\d w_c\over 2\pi i}
\prod_{1\le i< j\le c} \Delta^2(w_i)
{\prod_{\ell=1}^c \prod_{k=1}^{a+b}(w_\ell-{\gamma_k-\alpha\over\gamma_k-\beta})
\over \prod_{\ell=1}^c w_\ell^{b+c}}\cr
&={C\over c!}
\oint\cdots\oint {\d w_1\over 2\pi i}\ldots{\d w_c\over 2\pi i}
\prod_{1\le i< j\le c}\Delta(w_i)\Delta(w_i^{-1})
{\prod_{\ell=1}^c \prod_{k=1}^{a+b}(w_\ell-{\gamma_k-\alpha\over\gamma_k-\beta})
\over \prod_{\ell=1}^c w_\ell^b}
&\integrepc}$$
where $C=(\alpha-\beta)^{ab}\prod_{k=1}^{a+b}(\gamma_k-\beta)^c$.
The $w_i$ are integrated on contours surrounding $0$, for example $|w_i|=1$.
We recognize in the second line of Eq.~\integrepc\ the usual form of a matrix integral over the unitary group $U(c)$ once
angular variables are integrated out and only the eigenvalues $w_i$ are left. We thus obtain
\eqn\unitint{
\Phi_{a,b,c}=C' \int_{U(c)} \d\Omega\, \det(1+\Gamma\otimes\Omega) (\det \Omega^{-1})^{b}
}
where $\d\Omega$ is the Haar measure on $U(c)$, $\Gamma$ is the $(a+b)\times(a+b)$ diagonal matrix with eigenvalues 
${\gamma_k-\beta\over\alpha-\gamma_k}$,
and $C'=(\alpha-\beta)^{ab}\prod_{k=1}^{a+b}(\gamma_k-\alpha)^c$.
One can now use the identity $\det(1+\Gamma\otimes \Omega)=\sum_\lambda s_\lambda(\Omega) s_{\lambda^T}(\Gamma)$,
where $\lambda$ is a partition or Young diagram, $s_\lambda$ is the corresponding $GL$ character
(with the convention that it is zero if the Young diagram has more rows than the size of the matrix), 
and $\lambda^T$ is the transposed
Young diagram. Orthogonality of characters 
of $U(c)$ finally allows to perform the integration over $\Omega$
and results in the simple formula
\eqn\schurfn{
\Phi_{a,b,c}=C'\,s_{Y_{b,c}}(\Gamma)
}
where $Y_{b,c}$ is the rectangular Young diagram with $b$ rows and $c$ columns.
Of course, this is a standard result once reinterpreted in terms of non-intersecting paths 
($s_{Y_{b,c}}(\Gamma)$, as a Schur function of the eigenvalues of $\Gamma$,
can be defined as a sum over semi-standard Young tableaux,
which are themselves in bijection with the NIPs).

We finally consider the homogeneous
situation where all $z_i$ are equal to $1$, 
that is $\alpha_i=\alpha=1$, $\beta_j=\beta=q^2$, $\gamma_k=\gamma=q$. In this case
$\Gamma=-q\, 1$. We find that
$\Phi_{a,b,c}$ becomes $3^{(ab+bc+ca)/2}$ times 
the dimension
of the $GL(a+b)$ representation with rectangular Young diagram $b\times c$. The latter is 
one of the many formulae for the number of plane partitions.

The Razumov--Stroganov conjecture \RS\ claims that in the homogeneous $O(1)$
loop model,
$\Psi_{a,b,c}/\Psi_{n,min}$ (where $\Psi_{n,min}$ is the smallest component of $\Psi_n$) must be equal to the number of Fully Packed
Loop configurations (FPL) with the corresponding connectivity $(a,b,c)$. With our normalization conventions,
$\Psi_{n,min}=3^{n(n-1)/2}$
and all powers of $3$ cancel out, so that $\Psi_{a,b,c}/\Psi_{n,min}$ 
is simply the number of plane partitions in the hexagon $a\times b\times c$.
But according to \DFZJZ, the number of FPLs with connectivity $(a,b,c)$ is the very same number. This proves
the RS conjecture for the case of these link patterns.

\newsec{Generalization to four little arches}
Let us first reobtain the result of the previous section in a more synthetic
way. We use the fact, proved in appendix A, that 
when one switches the two spectral parameters of
neighboring parallel lines, the partition function of plane partitions with an arbitrary geometry
is unchanged.
This implies that the partition function $Z_{a,b,c}$ for plane
partitions introduced in Sect.~3 is a symmetric function of the spectral parameters $\{\alpha_i\}$, 
$\{\beta_j\}$, $\{\gamma_k\}$. At this stage,
one can skip the entire reinterpretation in terms of free fermions and
prove directly that it satisfies the same recurrence relations as $\Phi_{a,b,c}$. Indeed, setting
say $\alpha_1=\gamma_1$ forbids the lozenge parallel to sides $a$ and $c$ in the corner, see Fig.~\ppfrozen,
and thus creates two rows of ``frozen'' lozenges which lead us back to the case $a\times b\times (c-1)$.
This provides a nice graphical interpretation of the recurrence relations, very much in the spirit
of the recurrence relations of Korepin for the six-vertex model with Domain Wall Boundary Conditions \KOR.
\fig\ppfrozen{Plane partition in which the $\alpha_1$ and $\gamma_1$ rows are frozen.}{\epsfbox{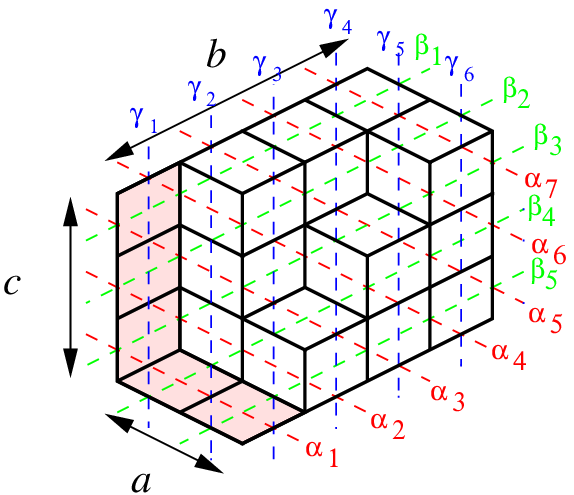}}

\fig\cinqpaq{Link pattern with four little arches $(a,b|e|c,d)$.}{\epsfxsize=5cm\epsfbox{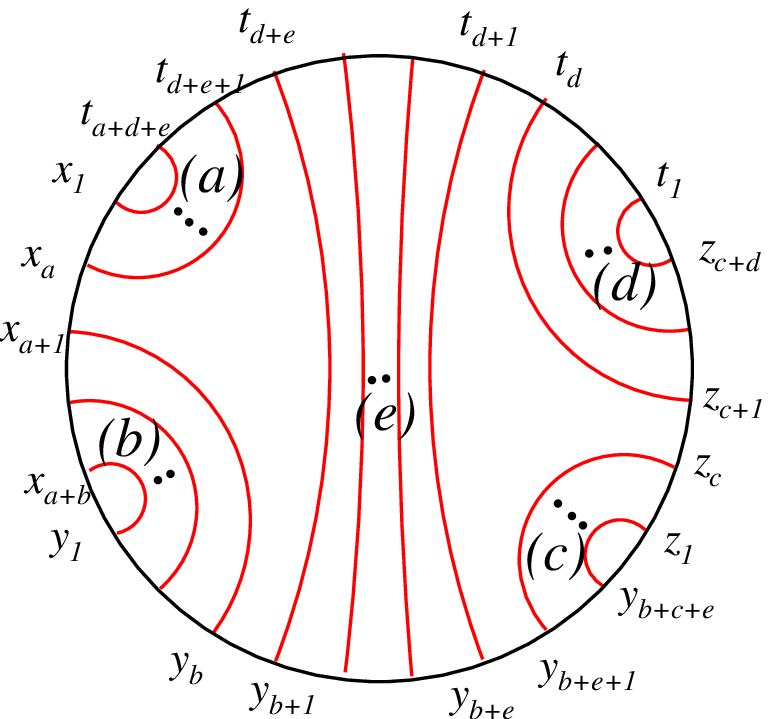}}
Consider now the most general link pattern which possesses four ``little arches'' (arches connecting neighbors),
as described by Fig.~\cinqpaq. It was shown in \DFZJZb\ that their enumeration is equivalent to that of certain
lozenge tilings of a region of the plane with identifications, see Fig.~\abcde. 
We refer the reader to \refs{\DFZ,\DFZJZb} for
details. In particular note that there are exactly $d$ ``dents'' in the two identified sides of length $c+d$.
Inspired by the case of three little arches, it is natural to introduce spectral parameters into
the lozenge tilings as described on Fig.~\abcde. The weight of a lozenge is equal to $q\, u - q^{-1}v$
where $u$ and $v$ are the spectral parameters crossing at the center of the lozenge in such a way
that the line of $v$ forms an angle of $+\pi/3$ with that of $u$ (contrary to the case of three little arches,
we cannot get rid of the factors of $q$ by a redefinition of the spectral parameters).
One checks that this produces a partition function
which has degree at most: $c+d+e$ in each $x_i$, $a+d$ in each $y_i$, $a+b+e$ in each $z_i$, $b+c$ in each $t_i$,
as should be. As a consequence of Appendix A, it is symmetric in each set of variables.
\fig\abcde{Lozenge tiling corresponding to the link pattern $(3,4|2|2,1)$
and its parameterization.}{\epsfbox{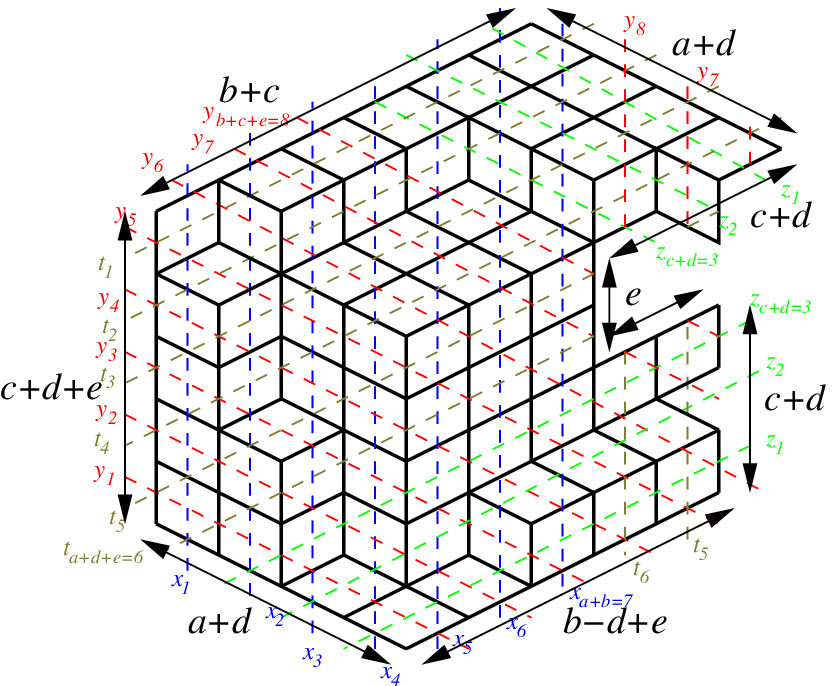}}

It is now easy to check that these partition functions satisfy all the required recurrence relations.
Among the various possibilities, two are depicted on Fig.~\abcdefr.
If one sets $t_1=q^2 x_1$, the first rows of the sides of length $a+d$ and $c+d+e$ are frozen,
and once these are removed one obtains the tiling with $a\to a-1$.
Similarly, if $t_1=q^{-2} z_1$, the first rows of the sides of lengths $b+c$ and $a+d$ are frozen, and
one dent becomes locked in first position. Once these rows are removed one recovers the tiling with
$d\to d-1$ (with, in particular, one less dent).
\fig\abcdefr{Lozenge tilings with $(t_1,x_1)$ and $(t_1,z_1)$ rows frozen.}{\epsfxsize=\hsize\epsfbox{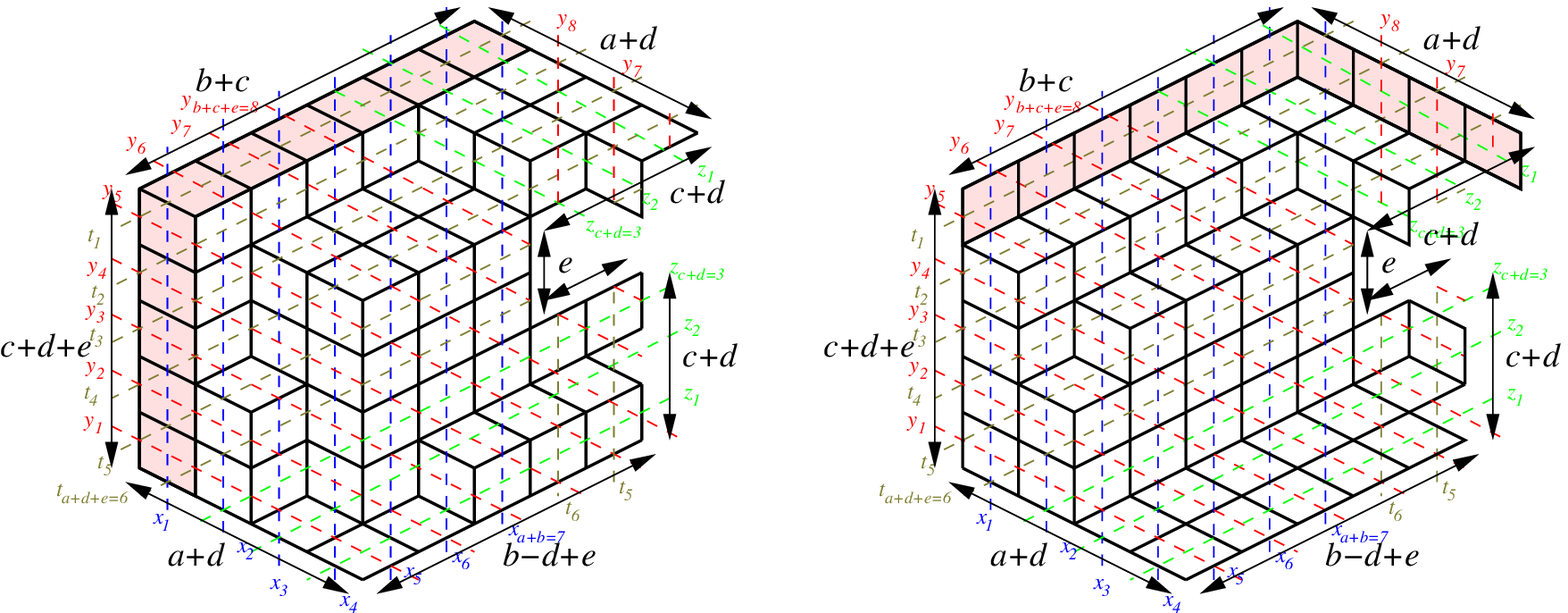}}
Since the partition function is of degree $b+c$ in $t_1$ and is known at $a+b+c+d$ values of $t_1$, it is entirely
fixed. 
The recurrence allows us to reach either $d=0$ or $a=0$, at which point, up to an additional frozen region
one recovers the tiling of a hexagon, a case which has already been treated.
Now the component of $\Psi_n$ corresponding to the link pattern $(a,b|e|c,d)$, once rid of its factors
$\prod_{i<j}(q\, x_i-q^{-1}x_j)\prod_{i<j}(q\, y_i-q^{-1}y_j)\prod_{i<j}(q\, z_i-q^{-1}z_j)\prod_{i<j}(q\, t_i-q^{-1}t_j)$,
satisfies the very same recurrence relations (of course at each step one must check
that the prefactors in the recurrence match); therefore it
is equal to the partition function of lozenge tilings described above. In particular, when all spectral
parameters are equal to $1$, the homogeneous component is equal to the number of such tilings, as predicted
by the RS conjecture.

\newsec{Conclusion}
In this article, we have proved in detail how the component $\Psi_{a,b,c}$ 
of the ground state of the $O(1)$ loop model
corresponding to the link pattern with three series of nested arches $(a,b,c)$,
once properly normalized, is equal to the number of plane partitions in a hexagon of size $a\times b\times c$.
This is a highly non-trivial check of the Razumov--Stroganov conjecture since these link patterns form an
infinite series. We have briefly described how the proof can be extended to more general link patterns,
$(a,b|e|c,d)$.

In fact, we have found more than this: just as in \DFZJ\ the sum of components 
was actually computed for arbitrary spectral parameters, here we have found
that $\Psi_{a,b,c}$ with spectral parameters is the partition function
of plane partitions with some local Boltzmann weights. This might seem
unsurprising
since one can hope that there is a unique ``natural'' way to introduce spectral
parameters into the model; however, if ones tries to connect to FPL 
configurations, in which the RS conjecture is formulated, one is faced
with a subtle problem: how to map back the plane partitions onto the square
lattice of the FPL in such a way as to make sense of the spectral parameter
dependence? Answering this is probably related to proving the full RS
conjecture. We hope to come back to this point in the future.

More obvious extensions of this work should be mentioned.
Firstly, more general FPL configurations have been \Kratt\ or could in principle be enumerated: they have
link patterns that can be obtained from those with four little arches by simple local modifications.
Clearly, the formulae of \DFZJ\ allow to obtain mechanically the corresponding ground state element.
It would be a further check of the RS conjecture to take the homogeneous limit in these formulae and
recover the counting of FPL, though it is unlikely to produce any surprise.

Secondly, An extension to arbitrary $q$ of the polynomials $\Psi_\pi$
was proposed in \Pas\ and reformulated as a solution of the $q$KZ equation
in \DFZJc. The present work being entirely based on recurrence relations
which are still satisfied by solutions of the $q$KZ equation, it is clear
that it can be extended to generic $q$. This might shed some light
on a possible generalization of the RS conjecture to arbitrary $q$.

\bigskip

\centerline{\bf Acknowledgments}
The author would like to thank J.-B. Zuber for carefully rereading the manuscript.

\bigskip

\appendix{A}{Proof of symmetry under interchange of spectral parameters}
The purpose of this appendix is to provide an explicit proof of the symmetry of the partition
function of weighted plane partitions when one switches the spectral parameters of two neighboring
parallel lines, say $\gamma_1$ and $\gamma_2$. The weighting is the same as described in Sec.~3 (i.e.\ for
simplicity we absorb powers of $q$ in the definition of the spectral parameters).
The plane partitions are supposed to fill a domain with arbitrary shape; in fact, we shall only
need to consider the region in which the two spectral parameters act and ignore the rest of the tiling.

Generally, this region has the shape of Fig.~\rtr. Note that the positions of the triangular holes are in
principle not fixed: to a given shape of the entire domain to be filled may correspond several
possibilities of locations for these holes. However
here we shall keep them fixed and show that each term in the summation over them is a symmetric function
of $\gamma_1$ and $\gamma_2$. The latter is equivalent to showing that a certain family of row-to-row transfer matrices is commutative; however here, to avoid any formalism we shall prove this by explicit computation.
\fig\rtr{Region of a tiling where the two spectral parameters $\gamma_1$, $\gamma_2$ appear.}{\epsfbox{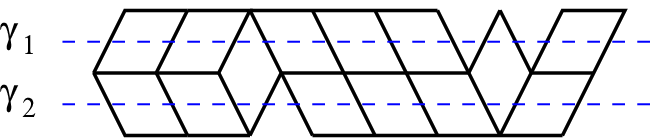}}

The locations of the triangular holes are strongly constrained in order to allow the possibility of a lozenge tiling. As a consequence we see that we can divide our region into pieces separated by holes, so that
these pieces have no
influence on each other, and the summation over plane partitions factorizes accordingly. There are two types
of pieces:

\item{a.}Hexagons of size $k\times 1\times 1$, $k\ge 0$. They are those that are delimited either (i) on at least one side by a pair of adjacent holes (by parity on the other side
the boundary must have the same shape) or (ii) by two individual holes, with the parity being such that
there are forced edges going inwards on both sides. The weighted enumeration of plane partitions inside a hexagon of size $k\times 1\times 1$ can be computed, cf Sec.~3:
\eqn\rega{
{\prod_{i=1}^{k+1} (\alpha_i-\gamma_1)(\beta_i-\gamma_2)-\prod_{i=1}^{k+1} (\alpha_i-\gamma_2)(\beta_i-\gamma_1)
\over \gamma_1-\gamma_2}
}
where the $\alpha_i$ and $\beta_i$ are the spectral parameters in the other 2 directions.
Eq.~\rega\ is explicitly symmetric in $\gamma_1$, $\gamma_2$.

\item{b.} Parallelograms $2\times k$, $k\ge 1$:
they are those that are delimited by two individual holes on opposite boundaries, with the
parity being such that there are forced edges going outwards. A parallelogram can only be filled by elementary
lozenges in a unique way, hence the weight
\eqn\regb{
\prod_{i=1}^k (\alpha_i-\gamma_1)(\alpha_i-\gamma_2)
}
where the $\alpha_i$ are the spectral parameters going parallel to the non-horizontal sides of the
parallelogram.

We finally conclude that the whole partition function, being a sum (over the locations of holes) of $\gamma_1$, $\gamma_2$ independent functions
(the partition function outside the region considered above) times products of functions of the type \rega\ or
\regb, is symmetric in $\gamma_1$, $\gamma_2$.

\listrefs
\end